\documentclass[12pt]{amsart}
\usepackage{amscd}

\title{Generalized Complete Intersections 
with Linear Resolutions}
\author{Munetaka Okudaira  and Yukihide Takayama}

\address{Munetaka Okudaira, Graduate School of Mathematical
Sciences, Ritsumeikan University, 
1-1-1 Nojihigashi, Kusatsu, Shiga 525-8577, Japan}
\email{rp001012@se.ritsumei.ac.jp}

\address{Yukihide Takayama, Department of Mathematical
Sciences, Ritsumeikan University, 
1-1-1 Nojihigashi, Kusatsu, Shiga 525-8577, Japan}
\email{takayama@se.ritsumei.ac.jp}
\theoremstyle{plain}
\newtheorem{thm}{Theorem}[section]
\newtheorem{prop}[thm]{Proposition}

\newtheorem{cor}{Corollary}
\theoremstyle{definition}

\theoremstyle{remark}

\newcommand{\G}[1]{{\rm G}(#1)}
\newcommand{\supp}[1]{{\rm supp}(#1)}
\newcommand{\F}{{\mathcal F}}
\newcommand{\mm}{{\mathfrak m}}

\let\Dirsum=\bigoplus
\let\iso=\cong
\def\NZQ{\Bbb}               % the font for N,Z,Q,R,C

\def\ZZ{{\NZQ Z}}

\begin{document}
\maketitle

\begin{abstract}
We determine the simplicial complexes $\Delta$ whose Stanley-Reisner 
ideals $I_{\Delta}$ have the following property: for all $n\geq 1$
the powers $I_{\Delta}^n$ have linear resolutions and finite length
local cohomologies.
MSC Primary: 13F55, Secondary: 13D02, 13D45
\end{abstract}

\section*{Introduction}

Let $S= K[X_1,\ldots, X_n]$ be a standard graded polynomial ring over a
field $K$ and $\mm = (X_1,\ldots, X_n)$. Recall that a graded ideal
$I\subset S$ is called {\em generalized Cohen-Macaulay} or simply {\em
FLC} (finite local cohomology) when the local cohomology
$H_{\mm}^{i}(S/I)$ has finite length for all $i< \dim S/I$.  We find
many examples of such ideals in algebraic geometry: the defining ideals
of Cohen-Macaulay projective schemes over the field are all FLC
ideals. However, as far as the authors are concerned,
we do not know very much about FLC monomial ideals.
The second author gave a combinatorial
characterization of FLC monomial ideals \cite{T} as an extension of the
well known case of squarefree monomial ideals, i.e., Buchsbaum
Stanley-Reisner ideals. 
On the other hand, the notion of generalized complete intersections (gCI) has 
been introduced in \cite{GT}. A gCI $I\subset S$ is a squarefree
monomial ideal such that the Stanley-Reisner ring $S/I$ is 
complete intersection over the punctured spectrum. 
All the powers $I^n$ ($n\geq 1$) of a gCI $I\subset S$ 
are FLC, and a combinatorical characterization of gCI has been given.

The purpose of this paper is to give a combinatorial characterization 
of the simplicial complexes $\Delta$ whose Stanley-Reisner ideals 
$I_\Delta$ are gCI and all the powers $I_\Delta^n$ ($n\geq 1$)
have linear resolutions. 
We will show that $I_\Delta$ is gCI and 
$I_\Delta^n$ has a linear resolution for all $n\geq 1$ if and only if 
$\Delta$ is as follows: $\Delta$ is a finite set of points,
the disjoint union of  paths or, if $\dim\Delta \geq 2$,
then $\Delta$ is a pure simplicial complex 
and is the disjoint union of facets $H$ and pairs $(F, G)$ of facets such
that $\vert F\setminus G\vert\; (=\vert G\setminus F\vert) =1$
(Theorem~\ref{thm:main}).
We also show that if the Buchsbaum
Stanley-Reisner ring $S/I_\Delta$  has minimal multiplicity 
in the sense of Goto \cite{G}, then any power $I_\Delta^\ell$, $\ell\geq 1$ 
is FLC and has a linear resolution (Corollary~\ref{thm:minmult}).
Finally, we give a computation of part of local cohomologies 
of gCI with linear resolutions (Proposition~\ref{thm:localcohomology}).

\section{Generalized Complete Intersection}

In this section, we recall some definitions and already known 
results, which will be used in the next section.

A Stanley-Reisner ring $K[\Delta] = S/I_\Delta$ is called
{\em generalized complete intersection} (gCI) if $K[\Delta]_P$ is 
a complete intersection for every prime ideal $P(\ne\mm)$ and $\Delta$ 
is a pure simplicial complex. For a simplicial complex $\Delta$, we will
always assume that $\{i\}\in \Delta$ for all $i\in [n]$.

\begin{thm}[cf. Th.~2.5 \cite{GT}]
Let $\Delta$ be a simplicial complex on the vertex set $[n] = \{1,\ldots, n\}$.
Then the following conditions are equivalent:
\begin{enumerate}
\item $K[\Delta]$ is a gCI;
\item $S/I_{\Delta}^{\ell+1}$ has FLC for an arbitrary integer $\ell \geq 0$.
\end{enumerate}
If one of these conditions holds, $K[\Delta]$ is Buchsbaum.
\end{thm}

A special case of gCI is complete intersection. In order to exclude this
uninteresting case, we use the notion of core of simplicial complex.
Let $\Delta$ be a simplicial complex on the vertex set $[n]$.
For $F \in \Delta$, we define ${\rm st}_{\Delta}(F) =
\{ G \in \Delta \;\vert\; G \cup F \in \Delta \}$.  We also define ${\rm
core}[n] = \{ i \in [n] \;\vert\; {\rm st}_{\Delta} (\{i\}) \ne \Delta
\}$. Then the {\it core} of $\Delta$ is defined by ${\rm core}\Delta =
\{ F \cap {\rm core}[n] \;\vert\; F \in \Delta \}$.
A gCI is a complete intersection if and only if 
${\rm core}\Delta\ne \Delta$. 

A combinatorial characterization of gCI is given as follows:
\begin{thm}[cf. Th.~ 3.16 \cite{GT}]
\label{thm:Goto-Takayama}
Let $K[\Delta]$ be a Stanley-Reisner ring with $\Delta = {\rm core}\Delta$
(i.e., $K[\Delta]]$ is not a complete intersection).
Let $\G{I_{\Delta}} = \{ u_1, \dotsc, u_{\ell} \}$ be the minimal
set of generators of $I_\Delta$ 
and $\F_{\Delta} = \{ \supp{u_j} \;\vert\; j = 1, \dotsc, \ell \}$, 
where $\supp{u} := \{j\in[n]\;\vert\; X_j\mbox{ divides }u\}$.
Then $K[\Delta]$ is a gCI
if and only if the following conditions hold:
\begin{enumerate}
\item $\Delta$ is pure.
\item for every $U\in\F_{\Delta}$ with $\vert U\vert\geq 3$, there
      exists a non-empty subset ${\mathcal C}(U)$ of $[n]$ such that 
  \begin{enumerate}
   \item ${\mathcal C}(U)\cap U=\emptyset$,
   \item for every $i\in {\mathcal C}(U)$, we have $E_{ij}:= \{i,j\}\in \F_{\Delta}$
   for all $j\in U$. Moreover if $U\cap T\ne\emptyset$ for
	 $T\in\F_{\Delta}$,
   then $T= E_{ij}$ for some $i,j$,
   \item for every $k\notin{\mathcal C}(U)\cup U$, we have
	 $\{i,k\}\in\F_{\Delta}$ for all $i\in{\mathcal C}(U)$.
  \end{enumerate}
\item Any two elements $i, j \in [n]$ are linked with a path 
$P = \{ \{i_k, i_{k+1} \} \;\vert\; k = 1, \dotsc, r \}$,
with edges $\{ i_k, i_{k+1} \} \in \F_{\Delta}$ for $k = 1, \dotsc, r$ 
such that $i = i_1$ and $j = i_{r+1}$.
\item If there exists a length $4$ path 
$P = \{ \{ i_p, i_{p+1} \} \in \F_{\Delta} \;\vert\; p = 1, 2, 3, 4 \}$
(with $i_1\ne i_{5}$), then there must be an edge 
$\{i_1, i_q \} \in \F_{\Delta}$ with $q = 3, 4$ or $5$.
\end{enumerate}
\end{thm}

Recall that a graded ideal $I\subset S$ is said to
have a linear resolution if all entries in the matrices representing the
differentials in a graded minimal $S$-free resolution of $I$ are linear.
It is an interesting question which ideal $I$ has the property that 
all the powers $I^n$ ($n\geq 2$) have a linear resolution.

An immediate corollary to Th.~\ref{thm:Goto-Takayama} is 
\begin{cor}
\label{thm:2-linear}
Let $I_{\Delta}\subset S$ is a gCI with ${\rm core}\Delta=\Delta$.
If $I_{\Delta}$ has a linear resolution, $I_{\Delta}$ is generated 
in degree $2$. Namely, a gCI may have only a 2-linear resolution.
\end{cor}
\begin{proof}If $I_{\Delta}$ has a $q$-linear resolution with $q\geq 3$.
Then $\deg u_i=q(\geq 3)$ for $i=1,\ldots, \ell$. But the condition 
$(ii)$ of Th.~\ref{thm:Goto-Takayama} implies the existence of a degree  $2$ 
element in $\G{I_{\Delta}}$, a contradiction.
\end{proof}

For $2$-linear resolution, we have the following result by 
Herzog-Hibi-Zheng.

\begin{thm}[Th.~3.2 \cite{HHZ}]\label{HHZheng}
Let $I$ be a monomial ideal generated in degree $2$. 
Then the following conditions are equivalent:
\begin{enumerate}
\item $I$ has a linear resolution.
%\item $I$ has linear quotients;
\item $I^\ell$ has a linear resolution for every $\ell\geq 1$.
\end{enumerate}
\end{thm}

In particular, when we consider Stanely-Reisner ideals,
a monimial ideal $I$ generated in degree~$2$ can be described 
in terms of edge graph. Namely, for a finite graph $G$ we 
define the edge ideal 
$I_G = (X_iX_j \;\vert\; \{i,j\}\mbox{is a edge of }G)$. Any
Stanley-Reisner ideal generated in degree~$2$ is the edge ideal
of a graph, which we call the {\em edge graph}.
Notice that when $I_\Delta$ is generated in degree~$2$, 
$\F_\Delta$ in Th.~\ref{thm:Goto-Takayama} 
can be identified with 
(the edge set of) the edge graph $G_\Delta$ corresponding 
to the edge ideal $I_\Delta$.

For linear resolutions of such ideals, we have 

\begin{thm}[Fr\"{o}berg \cite{F}]
\label{thm:Froberg}
Let $G$ be a graph. Then the edge ideal $I_G$ has a linear resolution
if and only if the complementary graph $\overline{G}$ is chordal.
\end{thm}
Recall that the {\em complemental graph} $\overline{G}$ of $G$ 
is the graph whose vertex set is the same as that of $G$ and 
whose edges are the non-edges of $G$. A graph $G$ is called 
{\em chordal} if each cycle of length $>3$ has a chord.
Notice that the complementary graph $\overline{G}$ in 
Th.~\ref{thm:Froberg} is exactly the 1-skeleton 
$\Delta_1$ of the simplicial complex $\Delta$ corresponding 
to the edge ideal $I_G$.

As an immediate consequence, we have 
\begin{cor}\label{thm:abstract}
Let $K[\Delta]$ be a Stanley-Reisner ring with $\Delta = {\rm core}\Delta$.
Then $K[\Delta]$ is a gCI and $I_\Delta^\ell$ has a linear 
resolution for every $\ell\geq 1$
if and only if the following conditions hold:
\begin{enumerate}
\item $\Delta$ is pure and $I_\Delta$ is generated in degree~$2$.
\item The 1-skeleton $\Delta_1$ is chordal.
\item Any two elements $i, j \in [n]$ are linked with a path 
$P = \{ \{i_k, i_{k+1} \} \;\vert\; k = 1, \dotsc, r \}$,
with edges $\{ i_k, i_{k+1} \} \in G_{\Delta}$ for $k = 1, \dotsc, r$ 
such that $i = i_1$ and $j = i_{r+1}$.
\item If there exists a length $4$ path 
$P = \{ \{ i_p, i_{p+1} \} \in G_{\Delta} \;\vert\; p = 1, 2, 3, 4 \}$
(with $i_1\ne i_{5}$), then there must be an edge 
$\{i_1, i_q \} \in G_{\Delta}$ with $q = 3, 4$ or $5$.
\end{enumerate}
\end{cor}

In the next section, we will give a precise description of the
simplicial complexes satisfying the conditions in
Cor.~\ref{thm:abstract}.

\section{Generalized complete intersections with 2-linear resolutions}
In this section, we will give a precise description of the
simplicial complexes $\Delta$ such that the Stanley-Reisner rings 
$K[\Delta]$ are gCI and every power $I_\Delta^\ell$ ($\ell\geq 1$) 
has a linear resolution. As we showed in the previous section,
this means gCI with a 2-liear resolution.

For a graph $G$ over the vertex set $V$, we define ${\rm Simp}(G)$ to be
the set of all subsets $F$ of $V$ such that $F$ is the vertex set of a
subgraph $H$ of $G$ isomorphic to a complete graph. For a simplicial
complex $\Delta$, ${\rm Simp}(\Delta_1)$, where $\Delta_1$ denotes the
1-skeleton, is the simpicial complex obtained by filling all the
simplicial cycles in $\Delta$.

\begin{prop}[cf. Prop.~6.1.25 \cite{V}]
\label{thm:Villarreal}
Let $I_\Delta$ be a Stanley-Reisner ideal generated.
Then $I_\Delta$ is generated in degree $2$
if and only if $\Delta = {\rm Simp}(\Delta_{1})$.
\end{prop}

\begin{prop}\label{thm:keylem}
Let $\Delta$ be a pure simplicial complex 
on the vertex set $[n]$ with $\Delta = {\rm core}\Delta$. 
Assume that $K[\Delta]$ is a gCI and $I_{\Delta}$ has 
a $2$-linear resolution.
Then, for any two distinct facets  $F$ and $H$ such that $F \cap H \ne \emptyset$, 
there exist a unique element $\{i_{1}, i_{2}\}\in G_{\Delta}$, i.e.,
$X_{i_1}X_{i_2}\in I_\Delta$,
with $i_{1} \in F \setminus H$ and $i_{2} \in H \setminus F$.
\end{prop}
\begin{proof}
Assume that there exists a pair of distinct facets $F$ and $H$  with 
$F\cap H\ne\emptyset$ such that no edge $\{i_1, i_2\}$, 
where $i_i\in F\setminus H$ and $i_2\in H\setminus F$, is in
$G_\Delta$. Then the complete graph ${\mathcal K}$ over
the vertex set $F\cup H$ is contained in the 1-skelton $\Delta_1$
of $\Delta$. But since $I_\Delta$ is generated in degree $2$, 
we know that $F\cup H\in \Delta$ by Prop.~\ref{thm:Villarreal},
which contradicts the assumption that $F$ and $H$ are facets. 
Thus we have proved the existence of the pair $\{i_1,i_2\}\in G_\Delta$
with the required property.

Now we show the uniqueness of the pair. Let $F$ and $H$ be distinct 
facets with $F\cap H\ne\emptyset$. We may assume 
$\vert F \setminus H\vert = \vert H\setminus F\vert  \geq 2$
and there exists a pair $\{i_1, i_2\}\in G_\Delta$ with
$i_1\in F\setminus H$ and $i_2\in H\setminus F$.

First of all, suppose there exists $j_{2} (\ne i_2) \in H \setminus F$
such that $\{i_{1}, j_{2} \} \in G_{\Delta}$.
Take any $i\in F\cap H$. Then, since $\Delta = {\rm core}\Delta$, there
exists $j\in [n]\setminus (F\cup H)$ such that $\{i, j\} \in
G_{\Delta}$. By Cor.~\ref{thm:abstract}(iii), 
there exists a path in $G_{\Delta}$ connecting $j$ and $i_{2}$, Furthermore,
by using Cor.~\ref{thm:abstract} (iv), we can take the path to be 
an edge  $\{j, i_{2}\}\in G_{\Delta}$. Consequently, we obtain 
the length $4$ path 
$\{ \{i, j\}, \{j, i_{2}\}, \{i_{2}, i_{1}\}, \{i_{1}, j_{2}\} \}$ 
in $G_\Delta$, so that by Cor.~\ref{thm:abstract}
at least one of  $\{i, i_{2}\}, \{i, i_{1}\}, \{i, j_{2}\} \in G_{\Delta}$. 
But since $\{i, i_{2}\}, \{i, j_{2}\} \in H$ and 
$\{i, i_{1}\}\in F$,
this cannot happen. Thus such $j_2$ does not exist.

Next suppose that there exist $j_{1} \in F \setminus H$ 
and $j_{2} \in H \setminus F$ such that
$j_{1} \ne i_{1}$, $j_{2} \ne i_{2}$ and $\{j_{1}, j_{2}\} \in G_{\Delta}$.
By the similar discussion as above, we know that 
$\{i_{1}, j_{2}\}, \{i_{2}, j_{1}\} \not \in G_{\Delta}$.
Thus $\Delta_1$ contains the length $4$ cycle 
$\{ \{i_{1}, j_{2}\}, \{j_{2}, i_{2}\}, \{i_{2}, j_{1}\}, \{j_{1}, i_{1}\} \}$
without any chord. But, since $I_\Delta$ has a $2$-linear resolution,
this contradicts Fr\"oberg's condition Th.~\ref{thm:Froberg}.
\end{proof}

\begin{cor}\label{thm:keycor}
Let $\Delta$ be as in Prop.~\ref{thm:keylem}.
Then, for any distinct facets $F$ and $H$ such that $F \cap H \ne \emptyset$,
we have $\vert F\setminus H\vert = \vert H\setminus F \vert =1$.
\end{cor}
\begin{proof}
By Prop.~\ref{thm:keylem}, there exist $i_{1} \in F \setminus H$ and $i_{2} \in H \setminus F$ 
such that $\{i_{1}, i_{2}\} \in G_{\Delta}$. By the uniqueness of $i_1$
and $i_2$, the 1-skeleton of $W:= (F \cup H) \setminus \{i_{1}, i_{2}\}$
is a complete graph so that $W \in \Delta$ by
Prop.~\ref{thm:Villarreal}.  Since $\Delta$ is pure, we have $W\subset U$ 
for some facet such that $\vert U\vert = \vert F\vert$.
Thus $\vert W \vert = \vert F\cup H\vert - 2 =
\vert F\vert + \vert H\vert - \vert (F\cap H) \vert -2 \leq \vert F\vert$,
so that we have $1\leq \vert H\vert - \vert (F\cap H)\vert \leq 2$.
Assume that $\vert H\vert  - \vert(F \cap H)\vert = 2$, which means 
that $W$ is a facet. Since $F\cap W\ne\emptyset$ and $F\setminus W = \{i_1\}$,
we have $W \setminus F = \{j\}$ for some $(i_2\ne)j\in (H\setminus F)$. 
Then by Prop.~\ref{thm:keylem} we must have $\{i_1, j\}\in G_\Delta$,
which contradicts the uniqueness of $\{i_1,i_2\}\in G_\Delta$ for 
the pair $F$ and $H$. Thus we must have 
$\vert H\vert  - \vert(F \cap H)\vert = 1$ as required.
\end{proof}

Now we are ready to prove our main theorem.

\begin{thm}
\label{thm:main}
Let $K[\Delta]$ be a gCI with $\dim K[\Delta] = d+1$ and 
${\rm core}\Delta = \Delta$.
Then $I_\Delta^\ell$ has a linear resolution for all $\ell\geq 1$
if and only if $\Delta$ is a finite set of points,i.e. $d=0$,
or otherwise $\Delta$ is as follows:
\begin{description}
\item [case ($\dim \Delta =1$)]
$\Delta$ is the disjoint union of paths $\Gamma_1, \ldots, \Gamma_s$
of arbitrary lengths.
\item [case ($\dim \Delta \geq 2$)]
$\Delta$ is the disjoint union of the following types of 
subcomplexes:
\begin{description}
\item [type 1] $\langle F, G\rangle$, where $F$ and $G$ are $d$-simplexes
such that $F\cap G\ne\emptyset$ 
and $\vert F\setminus G\vert = \vert G\setminus F\vert = 1$,
\item [type 2] $\langle H\rangle$, where $H$ is a $d$-simplex.
\end{description}
Notice that, since ${\rm core}\Delta =\Delta$, we exclude the case 
that $\Delta$ itself is of type~1 or type~2.
\end{description}
\end{thm}
\begin{proof}
Notice that if $\dim \Delta =0$, then $K[\Delta]$ is clearly a
Cohen-Macaulay gCI and since $I_\Delta = (X_iX_j\;\vert\; 1\leq i<j<n)$
we know that $I_\Delta$ has a $2$-linear resolution so that by 
Th.~\ref{HHZheng} $I_\Delta^\ell$ also has a linear resolution
for all $\ell\geq 1$. In the following, we consider the case of $d>0$.
Notice that since $K[\Delta]$ is gCI, $\Delta$ is pure.
If $\Delta$ is as stated above, it is straightforward to check that 
$K[\Delta]$ satisfies the conditions of Cor.~\ref{thm:abstract}.

Now assume that $I_\Delta^\ell$ has a linear resolution for all
$\ell\geq 1$.  Then by Cor.~\ref{thm:2-linear} $I_\Delta$ is generated
in degree~$2$.  

Consider any distinct facets $F$ and $G$ with $F\cap G\ne\emptyset$ 
and assume that there exists the third facet $H$ with $H \cap (F\cup G)\ne\emptyset$.
By Cor.~\ref{thm:keycor}, we know that we must have 
$C:= F\cap G = G\cap H = F\cap H \ne \emptyset$
and $\{i_i\} = F\setminus C$, $\{i_2\} = G\setminus C$ and $\{i_3\} = H \setminus C$
for some distinct $i_1$, $i_2$ and $i_3$. Let $i\in C$ be arbitrary.
Then, there exists $j\notin F\cup G\cup H$ since ${\rm core}\Delta = \Delta$.
Now, as in the proof of Prop.~\ref{thm:keylem}, we obtain the length $4$ path
$\{i,j\}, \{j, i_2\}, \{i_2, i_3\}, \{i_3, i_1\}$ in 
the edge graph $G_\Delta$, for which $\{i, i_1\}, \{i, i_2\}, \{i, i_3\}\notin G_{\Delta}$
since each of them is in a facet. This contradicts
Cor.~\ref{thm:abstract} (iv).  Thus such a facet $H$ does not 
exist. This implies that $\Delta$ is as stated above.
\end{proof}

\begin{cor}
Let $\Delta$ be a pure simplicial complex with ${\rm core}\Delta = \Delta$. 
Then the power of the Stanely-Reisner ideal $I_\Delta^\ell$ 
has FLC and a linear resolution for 
all $\ell\geq 1$, if and only if $\Delta$ is 
as in Th.~\ref{thm:main}.
\end{cor}

\section{Stanley-Reisner ring with minimal multiplicity}

Let $A = K[A_{1}]$ be a homogeneous Buchsbaum $K$-algebra of 
dimension $d$ with the unique homogeneous maximal ideal
$\mm = A_{+}$. Then $A$ is called a {\it Buchsbaum ring with minimal
multiplicity} (\cite{G}) if 
$e(A) = 1 + \sum_{i=1}^{d-1} \binom{d-1}{i-1}l_{A}( H_{\mm}^{i}(A))$, 
where $e(A)$ denotes
the multiplicity of $A$ with regard to $\mm$,
$H_{\mm}^i(A)$ is the $i$th local cohomology and 
$\ell_{A}(M)$ denotes the length of the $A$-module.

On the other hand, a Stanley-Reisner ideal $I_\Delta \subset S$ 
generated in the same degree $\delta$ is called {\em matroidal} 
of degree $\delta$ if 
$\mathcal{B}:= \{\supp{u}\;\vert\; u\in \G{I_\Delta}\}$
forms a matroid base of degree $\delta$. Namaly,
$\vert B\vert = \delta$ for all $B\in\mathcal{B}$, and 
for all $B_1, B_2 \in \mathcal{B}$ and 
all $i\in B_1 \setminus B_2$, 
there exists 
$j\in B_2 \setminus B_1$ such that
$(B_1 \setminus \{i\}) \cup \{j\} \in \mathcal{B}$.
In particular, for a Stanley-Reisner ideal $I_\Delta$ generated 
in degree~$2$ is matroidal if the edge graph $G$ has the following 
property: for any disjoint edges $\{i,j\}$  and $\{p,q\}$ 
of the edge graph $G_\Delta$,
each vertex $\{k\}$, $k = i, j$ is linked with at least 
one of the vertices of the edge $\{p,q\}$, and vice versa.
Matroidal ideals have linear resolutions.
See \cite{HT} for the detail of this fact. 

We note that, for a matroidal ideal $I_\Delta$, the corresponding
simplicial complex is not always pure. For example, 
\begin{equation*} 
I_{\Delta} = (X_{1}X_{3}, X_{1}X_{4}, X_{2}X_{3}, X_{2}X_{4},
X_{3}X_{5}, X_{4}X_{5}),
\end{equation*}
is matroidal, and the simplicial complex $\Delta$ is spanned by
the facets $\{1,2,5\}$ and $\{3,4\}$.

Now we show the following.

\begin{prop}
Let $K[\Delta]$ be a Buchsbaum Stanely-Reisner ring with $\dim K[\Delta] =d+1$
and ${\rm core}\Delta = \Delta$.
Then the followings are equivalent:
\begin{enumerate}
\item [$(i)$] $K[\Delta]$ has minimal multiplicity.
\item [$(ii)$] $\Delta$ is the disjoint union of 
$d$-simpliexes.
\item [$(iii)$] $I_\Delta$ is matroidal of degree~$2$.
\end{enumerate}
\end{prop}
\begin{proof}
The equivalence of $(i)$ and $(ii)$ is shown in 
Example~3.1 \cite{TY}. 
Assume that $\Delta = \langle H_1,\ldots, H_r\rangle$
(disjoint union of $d$-simplexes).
Then $I_\Delta = (X_iX_j\;\vert\; i\in H_p,\; j\in H_q
\mbox{ for some }1\leq p < q \leq r)$, so that we easily
know that this ideal is matroidal. Now we have only to
show $(iii)$ to $(ii)$. 

$\Delta$ is pure  since $K[\Delta]$ is Buchsbaum.
Assume that $I_\Delta$ is matroidal
of degree~$2$ but $\Delta$ is not the disjoint union of $d$-simplexes.
Then there exist facets $F$ and $G$ such that $F\cap G\ne \emptyset$.
Choose any $i\in F\cap G$.
Since ${\rm core}\Delta=\Delta$, there exists $j\notin F\cup G$
such that $\{i, j\}\in G_\Delta$.

On the other hand, by Th.~\ref{thm:Goto-Takayama}, we easily know that 
$I_\Delta$ is a gCI.  Thus 
there exists 
$\{i_1,i_2\}\in G_\Delta$ with $i_1\in F\setminus G$ 
and $i_2\in G\setminus F$ by Prop.~\ref{thm:keylem}.

Since $\{i_1, i\}, \{i_2,i\}\notin G_\Delta$, the existence of 
two disjoint edges $\{i,j\}$ and $\{i_1, i_2\}$ contradicts 
the assumption that $I_\Delta$ is matroidal.
\end{proof}

\begin{cor}
\label{thm:minmult}
Let $K[\Delta]$ be a Buchsbaum Stanley-Reisner ring with
minimal multiplicity.
Then $I_\Delta^\ell$ has FLC and 
a linear resolution for for all $\ell\geq 1$.
\end{cor}

\section{Local cohomologies of generalized complete intersection
 with linear resolutions}

In this section, we consider local cohomologies of the generalized
complete intersection with linear resolutions.  For a positively graded
$K$-algebra $R$ with 
the graded maximal ideal
$\mm = \Dirsum_{n>0}R_n$, we denote by $H^i_\mm(R)$
the $i$th local cohomology module with regard to 
$\mm$.  Since we consider monomial ideals, the local cohomology
modules have the $\ZZ^n$-grading. In the following, we will denote
by $[H^i_\mm(R)]_a$, where $a\in\ZZ$ or $a\in\ZZ^n$, the $a$-th 
graded component of the module.

It is well known that, for Stanely-Reisner ideals, FLC and Buchsbaum are
equivalent notions and for a Buchsbaum Staneley-Reisner ring
$K[\Delta]$, we have $[H^i_\mm(K[\Delta])]_j=0$ for all $i<\dim
K[\Delta]$ and for all $j\ne 0$. For monimial ideals, if $I\subset S$
is FLC then $[H^{i}_\mm(S/I)]_j=0$ for all $i<\dim S/I$ and $j<0$
or $\sum_{k=1}^n\rho_k - n <j$ where $\rho_k$ is the maximal exponent
of the variable $X_k$ in the minimal set $G(I)$ of monomial generators.
In particular, we have 
\begin{prop}[Prop.~1 \cite{T}]
\label{thm:monoFLC}
For a monomial ideal $I\subset S$, the local cohomology $H^i_\mm(S/I)$,
$i\ne \dim S/I$, has finite length if and only if 
$H^i_\mm(S/I)_a=0$ for all $a\in\ZZ^n$ such that 
$a_j<0$ for some $1\leq j\leq n$.
\end{prop}
See \cite{T} for the results on FLC monomial ideals.

We will now compute $[H^i_\mm(S/I_\Delta^\ell)]_0$ for all $i<\dim
S/I_\Delta^\ell = \dim S/I_\Delta$ and all $\ell \geq 1$, for 
gCI $I_\Delta$  with a linear resolution.

We first recall a few results.
The local cohomologies of a monomial ideal $I$
and its radical $\sqrt{I}$ can be compared by
\begin{prop}[cf. Cor.~2.3 \cite{HTT}]
\label{thm:compare}
Let $I\subset S$ be a monimial ideal. Then we have 
the following isomorphisms of $K$-vector spaces
\begin{equation*}
   [H^i_\mm(S/I)]_a \iso [H^i_\mm(S/\sqrt{I})]_a
\end{equation*}
for all $a\in\ZZ^n$ with $a_i\leq 0$ for all $1\leq i\leq n$.
\end{prop}

For the local cohomologies of Stanely-Reisner ideal, we recall
Hochster's formula:
\begin{thm}[cf. Th.~5.3.8 \cite{BH}]
\label{thm:hochster}
The Hilbert series of $K[\Delta]$ with respect to 
the $\ZZ^n$-grading is given by
\begin{equation*}
{\rm Hilb}(H^i_\mm(K[\Delta]), {\bf t})
= \sum_{F\in\Delta}
   \dim_K \tilde{H}_{i-\vert F\vert -1}({\rm lk} F; K)
   \prod_{j\in F}\frac{t_j^{-1}}{1- t_j^{-1}}.
\end{equation*}
\end{thm}

Now we show the following

\begin{prop}
\label{thm:localcohomology}
Let $I_\Delta \subset S$ be a gCI with a linear resolution.
Then we have the following isomorphism as $K$-vector spaces
\begin{equation*}
[H^i_\mm(S/I_\Delta^\ell)]_0
\iso \left\{
      \begin{array}{ll}
       K^{\alpha -1} & \mbox{if $i =1$}\\
       0  & \mbox{if $i\ne 1,d$}\\
      \end{array}
      \right.
\end{equation*}
where $d = \dim S/I_{\Delta} = \dim \Delta +1$
and $\alpha$ is the number of connected components 
of $\Delta$.
\end{prop}
\begin{proof} 
By Prop.~\ref{thm:monoFLC} and Prop.~\ref{thm:compare},
we have $[H^i_\mm(S/I_\Delta^\ell))]_0 
= [H^i_\mm(S/I_\Delta^\ell))]_{{\bf 0}}
= [H^i_\mm(S/\sqrt{I_\Delta^\ell})]_{{\bf 0}}
= [H^i_\mm(S/I_\Delta)]_{{\bf 0}}$
for all $i<d$ and $\ell \geq 1$,
where ${\bf 0} = (0,\ldots, 0)\in\ZZ^n$.
Now, since $I_\Delta \subset S$ is a gCI with a $2$-linear resolution,
we compute $[H^i_\mm(S/I_\Delta)]_0=0$ for $i\ne 1$ and $i<d$ using
Th.~\ref{thm:hochster} (see Prop.~1.1(3) in \cite{TY}). Now we 
consider the case of $i=1$.
Since $\Delta$ is as in Th.~\ref{thm:main},
${\rm lk}(F)$ ($F\in\Delta$) and 
${\mathcal H}:= \tilde{H}_{i-\vert F\vert -1}({\rm lk}(F);K)
= \tilde{H}_{\vert F\vert}({\rm lk}(F);K)$
are as follows:
\begin{enumerate}
\item ${\rm lk}(F)=\Delta$ if $F=\emptyset$.
Then ${\mathcal H} \iso K^{\alpha-1}$, where $\alpha$ 
is the number of connected components in $\Delta$.
\item ${\rm lk}(F) = \{\emptyset\}$, if $F$ is an isolated facet.
Then ${\mathcal H} = 0$.
\item ${\rm lk}(F)$ is a simplex, if $F$ is one of the following;
\begin{enumerate}
\item $F$ is a non-facet of an isolated facet,
\item for a pair $\langle F', G'\rangle$ of intersecting facets
we have $F\subset F'$ and  also  $F'\setminus G' \subset F$.
\end{enumerate}
Then ${\mathcal H} = 0$.
\item ${\rm lk}(F)$ is two points, if $F$ is the intersection of two facets 
$F'$ and $G'$ such that $F'\cap G' = F$ and $\vert F'\setminus G'\vert
= \vert G'\setminus F'\vert =1$.
Then, ${\mathcal H} = 0$.
\item ${\rm lk}(F)$ is a union of two simplexes $\langle F', G'\rangle$,
where $F'\cup F$ and $G'\cup F$ are facets.
Then, ${\mathcal H} = 0$.
\end{enumerate}
Then we obtain the desired result by Th.~\ref{thm:hochster}.
\end{proof}

\end{document}